\documentclass[11pt, psamsfonts]{amsart}
\usepackage{amssymb, amsmath, amsthm}
\usepackage[final, hypertex]{hyperref}
\usepackage[a4paper, centering]{geometry}
\newtheorem{Theorem}{Theorem}[section]
\newtheorem{Proposition}[Theorem]{Proposition}
\newtheorem{Lemma}[Theorem]{Lemma}
\newtheorem{Corollary}[Theorem]{Corollary}
\theoremstyle{definition}
\newtheorem{Definition}[Theorem]{Definition}
\theoremstyle{remark}
\newtheorem{Remark}[Theorem]{Remark}

\long\def\elimina#1{} 

\def\R{\mathbb{R}}

\def\N{\mathbb{N}}
\def\U{\mathcal{U}}
\def\Wuu{W^{1,1}_0}
\def\inrad#1{r_{#1}}
\def\inradius{\inrad{\Omega}}
\def\pscal#1#2{\left\langle#1,\,#2\right\rangle}

\def\convopen{\mathcal{K}^n_0}

\def\gau#1{\rho_{#1}}
\def\gauge{\rho}

\def\pgauge{\gauge^0}

\def\dist{d}

\def\cut{m}

\def\curvg{\tilde{\kappa}}

\def\nor{\nu}

\def\vf{v_f}
\def\meas{\mathcal{L}^n}

\def\inrad#1{r_{#1}}
\def\inradius{\inrad{\Omega}}
\def\Wuu{W^{1,1}_0}

\def\mean{\widetilde{H}}

\def\pscal#1#2{\langle#1,\,#2\rangle}

\def\nor{\nu}
\def\Lip{\textrm{Lip}}
\def\Lipr{\Lip^1(\Omega,\rho)}
\def\Liprz{\Lip^1_0(\Omega,\rho)}
\def\TS{T_S}

\def\TSe{T_S^e}
\def\Tf{T_f}
\def\oray#1{{]\!]#1[\![}}
\def\cray#1{{[\![#1]\!]}}
\def\lray{\lambda^*}
\def\ch{n_{\gauge}}

\DeclareMathOperator{\inte}{int} \DeclareMathOperator{\spt}{supp}
\DeclareMathOperator{\conv}{co} \DeclareMathOperator{\dive}{div}
 
\DeclareMathOperator{\proj}{\Pi}

\begin{document}
\title[Variational problems solved by a distance function]
{A sharp uniqueness result for a class of variational problems solved by a distance function}%

\author[G.~Crasta]{Graziano Crasta}
\address{Dipartimento di Matematica ``G.\ Castelnuovo'', Univ.\ di Roma I\\
P.le A.\ Moro 2 -- 00185 Roma (Italy)}
\email[Graziano Crasta]{crasta@mat.uniroma1.it}

\author[A.~Malusa]{Annalisa Malusa}
\email[Annalisa Malusa]{malusa@mat.uniroma1.it}

\date{December 20, 2006}

\keywords{Minimum problems with constraints, uniqueness, Euler equation, distance function,
mass transfer problems, $p$-Laplace equation}
\subjclass[2000]{Primary 49J10}

\dedicatory{(Dedicated to Arrigo Cellina in the occasion of his 65-th birthday)}

\begin{abstract}
We consider the minimization problem for an integral functional
$J$, possibly non-convex and non-coercive in $W^{1,1}_0(\Omega)$,
where $\Omega\subset\R^n$ is a bounded smooth set.
We prove sufficient conditions in order to guarantee that
a suitable Minkowski distance is a minimizer of $J$.
The main result is a necessary and sufficient condition
in order to have the uniqueness of the minimizer.
We show some application
to the uniqueness of solution of a system of PDEs of
Monge-Kantorovich type arising in problems of mass transfer theory.
\end{abstract}

\maketitle

\section{Introduction}

In this paper we consider the minimization problem
\begin{equation}\label{f:J}
\min_{u\in\Wuu(\Omega)}J(u),
\qquad J(u) = \int_{\Omega} [h(Du) - f\, u]\,
dx\,,
\end{equation}
where $\Omega\subset\R^n$ is a bounded smooth (i.e.\ $C^2$) open
set, $h\colon\R^n \to [0,+\infty]$ is a (possibly non-convex)
Borel function, and $f\in
L^{\infty}(\Omega)$ is a non-negative function.
We are interested in existence and uniqueness results.

The main assumptions on $h$ involve
the convex hull $K$ of its zero-level set
\[
K :=\conv Z,\qquad
Z := \{\xi\in\R^n;\ h(\xi)=0\},
\]
and the rate of growth of $h$ outside $K$,
defined by
\begin{equation}\label{lak}
\Lambda := \sup \left\{\lambda\geq 0;\
h(\xi)\geq \lambda (\rho (\xi)-1)\ \forall \xi\in \R^n\right\}\,,
\end{equation}
where $\gauge\colon\R^n\to\R$ is the gauge function
of the convex set $K$ (see Section \ref{ss:conv}).
More precisely, the assumptions on $h$ are the following:
\begin{itemize}
\item[(H1)] $h\colon \R^N \to [0,+\infty]$ is a Borel function;

\item[(H2)]
$K$ is a compact convex set containing $0$ as an interior point,
with boundary $\partial K$ of class $C^2$, and strictly positive principal curvatures.
\end{itemize}

Notice that (H2) implies that $K$ is a strictly convex set. Hence $\partial K\subseteq Z$,
and $h(\xi) = 0 = \min h$ for every $\xi\in\partial K$.

The main tool needed in our investigation is
the Minkowski distance $\dist$ from $\partial\Omega$ associated to $K$
(see Section \ref{ss:dist}).
It is well known that $\dist$ is a Lipschitz continuous
function vanishing on $\partial\Omega$,
with $D\dist \in \partial K$ a.e.\ in $\Omega$.
In fact, $\dist$
is the maximal element
of the family $\Liprz$ of
all Lipschitz functions
$u\colon\overline{\Omega}\to\R$
satisfying $Du\in K$ a.e.\ in $\Omega$ and $u=0$ on $\partial\Omega$.
Since $f\geq 0$, a direct inspection of $J$
suggests that it is reasonable to expect
that $\dist$
is a minimizer of $J$, provided that $h$ grows fast enough
outside $K$.
This guess will be proved in Theorem~\ref{t:exif} below.

The aim of this paper is to give necessary and sufficient
conditions in order to guarantee that $\dist$ is the unique minimizer
of $J$. The leading idea of our techniques is to
go deeply into the connection between these problems
of Calculus of Variations and some aspects of the mass
transport theory.
In order to explain this connection,
let us consider the model functional
\begin{equation}\label{f:celp}
J(u)= \int_{\Omega} [h(|Du|) - u]\,
dx\,,
\end{equation}
studied by A.~Cellina in the seminal paper \cite{Ce1}.
Here $\Omega$ is a \textsl{convex} subset of $\R^2$, and the Lagrangian
is radially symmetric, so that $K$ is a ball.
Assume, without loss of generality,
that $K=\overline{B}_1(0)$. In this case $\rho(\xi)=|\xi|$,
and $\dist$ is the Euclidean distance function from the
boundary of $\Omega$.
In \cite{Ce1} it was proved that $\dist$ is a solution
to (\ref{f:celp}), provided that the growth condition $\Lambda\geq\inradius$
holds (here $\inradius$ is the inradius of the set $\Omega$).
The key point in the proof of the minimality of $\dist$ is to
observe that for every $v\in L^\infty(\Omega)$, satisfying $0\leq v(x)\leq \Lambda$
a.e.\ in $\Omega$, one has
\begin{equation}\label{f:estintr}
\begin{split}
J(u)-J(\dist)&\geq \int_\Omega \left[v\, \max\{|Du|-1, 0\}-u+\dist\right]\, dx \\
& \geq \int_\Omega \left[v \pscal{D\dist}{Du-D\dist}-(u-\dist)\right]\, dx\,.
\end{split}
\end{equation}
Hence the conclusion follows once one exhibits a function $v$ as above which solves
\[
-\dive(v\, D\dist)=1\quad \textrm{in}\ \Omega\,,
\]
in the sense of distributions.
In fact, the hard part of the proof of the existence result
in \cite{Ce1} is the construction of such a function $v$.
Moreover it is shown, by examples, that
the growth condition $\Lambda\geq\inradius$
cannot be improved.
The result in \cite{Ce1} has been extended to convex
domains in $\R^n$ and to more general functionals in subsequent works (see
\cite{CeCe,CPT,CMd,Tr,Vo}).

Recently in \cite{CaCa,CCCG} it was proved that for every given non-negative continuous function
$f$ there exists a unique non-negative continuous function $\vf$ solving
\[
-\dive(\vf\, D\dist)=f\quad \textrm{in}\ \Omega\,,
\]
without the requirement of $\Omega$ to be a convex set.
In dealing with non-convex domains, the growth
condition has to be modified,
in order to take into account the presence of points on $\partial \Omega$
with negative curvatures.
Nevertheless it can be proved that,
if ${\|f\|}_{\infty}$ is small enough
(see (H3) below),
then the function $\vf$ satisfies $0\leq \vf <\Lambda$ and hence
$\dist$ is a solution to the minimum problem
\begin{equation}\label{f:celf}
\min_{u\in\Wuu(\Omega)}\int_{\Omega} [h(|Du|) - fu]\, dx\,.
\end{equation}
Going further in the study of Cellina's minimization problem, one easily get
that the estimates $0\leq \vf <\Lambda$ imply an \textit{a priori} bound on the
gradient of minimizers.
Namely, as a consequence of the analogous of (\ref{f:estintr}),
with $f$ instead of 1, every
solution $u\in\Wuu(\Omega)$ of (\ref{f:celf}) has to belong in fact
to $W^{1,\infty}_0(\Omega)$, and $|Du|\leq 1$.
Hence some information on the solutions to
(\ref{f:celf}) can be obtained by studying the ancillary minimization problem with
constraints
\begin{equation}\label{f:minfac}
\min_{u\in\Wuu(\Omega),\ |Du|\leq 1}-\int_{\Omega}fu\, dx\,,
\end{equation}
which fits into a branch of the optimal mass transfer theory.
It is plain that
$\dist$ is always a solution of (\ref{f:minfac}) for every bounded set and for
every $f\in L^1(\Omega)$, $f\geq 0$ a.e.\ in $\Omega$ (and even for non-negative
bounded measures).
Clearly, if the (essential) support $\spt (f)$ of $f$ coincides with ${\Omega}$,
then $\dist$ is the unique minimizer.
Moreover, it is well known that in the optimal mass transfer
problems a long--range effect occurs, i.e.\ $\dist$ is the unique solution to (\ref{f:minfac})
if the mass displacement spreads over the whole $\Omega$ (see \cite{CaCa,CCCG}).
We shall show that, in fact, $\dist$ is the unique solution to (\ref{f:minfac}) also in many
cases where the mass transfer spreads only a part of $\Omega$.

\smallskip
The arguments above remain valid also in the anisotropic case.
A key point in the analysis of problem (\ref{f:J})
is the study of the constrained minimization problem
\begin{equation}\label{f:minfak}
\min_{u\in\Wuu(\Omega),\ Du\in K}-\int_{\Omega}fu\, dx.
\end{equation}
We shall show that the Minkowski distance $d$ associated to $K$
is the unique solution to (\ref{f:minfak}) if and only if
$\spt(f)$ contains the singular set $\Sigma$ of those points
where $\dist$ is not differentiable.
Furthermore,
we exhibit an explicit solution $u_f$ to (\ref{f:minfak})
with $D u_f\in\partial K$ a.e.\ in $\Omega$,
that coincides with $\dist$
if and only if $\Sigma\subseteq\spt(f)$
(see Theorem~\ref{t:uf} below).
These results lead to the fact that problem (\ref{f:J})
has either $\dist$ as unique solution, if $\Sigma\subseteq\spt(f)$,
or it has at least two distinct solutions $u_f$ and $\dist$,
if $\Sigma\setminus\spt(f)\neq\emptyset$.

The r\"ole of the singular set $\Sigma$ in the uniqueness result
can be understood from a mass transfer theory viewpoint.
Namely, if $f$ is a non-negative continuous function, it can be proved that $u$ is a solution of the
constrained minimization problem (\ref{f:minfak})
if and only if there exists a non-negative continuous function $v$ such that the
pair $(u,v)$ is a solution to the system of PDEs
\begin{equation}\label{f:syst1}
\begin{cases}
-\dive(v\, D\gauge(Du)) = f
&\textrm{in $\Omega$},\quad \textrm{(distributional)}\\
\gauge(Du)\leq 1
&\textrm{in $\Omega$},\\
\gauge(Du) = 1
&\textrm{in $\{v>0\}$}\quad \textrm{(viscosity)},
\end{cases}
\end{equation}
complemented with the conditions
\begin{equation}\label{f:syst2}
\begin{cases}
u\geq 0,\
v\geq 0
&\textrm{in $\Omega$},\\
u=0
&\textrm{on $\partial\Omega$}
\end{cases}
\end{equation}
(see Section~\ref{s:MK}).
In \cite{CMf,CMg}
we have constructed a
non-negative continuous function $\vf$
such that $(\dist,\vf)$ is a solution of (\ref{f:syst1})--(\ref{f:syst2}).
Moreover we have proved that if $(u,v)$ solves (\ref{f:syst1})--(\ref{f:syst2}), then
$v= \vf$, and $u=\dist$ in $\{\vf>0\}$.
In terms of optimal transport problems, the conclusions
above state that, for every given mass density $f\geq 0$, the transport density $\vf$
is uniquely determined, while the transport potential $u$ may differ from $\dist$
only in the region $\{\vf=0\}$ where no mass transfer occurs (long range effect).

For what concerns the uniqueness of the solution,
we start from the fact that
$\dist$ is the unique element $u$ in the family $\Liprz$
matching the condition $u=\dist$ on $\Sigma$.

The results in \cite{CMg} imply that if $u$ is a solution of (\ref{f:J}),
then $u=\dist$ in the transport set $\{\vf>0\}$.
Furthermore, due to its structure, the closure of the transport set
contains $\Sigma$ if and only
if $\Sigma\subseteq\spt(f)$.
Then, whenever $\Sigma\subseteq\spt(f)$,
$\dist$ is the unique solution to (\ref{f:J}).

\smallskip
The paper is organized as follows.
In Section~\ref{s:prel}
we introduce some notation and review some preliminary
result.
In Section~\ref{s:exis} we recall the properties
of the transport density $\vf$, associated to the
Monge-Kantorovich system (\ref{f:syst1})-(\ref{f:syst2}),
and we prove that the Minkowski distance $\dist$
is a solution to the minimization problem (\ref{f:J}),
provided that ${\|f\|}_{\infty}$ is small enough.
Section~\ref{s:uni} is devoted to the proof of the
uniqueness result.
More precisely, in the first part of the section
we prove some fine property
of functions defined by a max-convolution
of cone-shaped functions.
These properties will be essential in order to obtain
necessary and sufficient conditions for the uniqueness
of the minimizer in problem~(\ref{f:J}).
Finally, in Section~\ref{s:ext}
we give some extension to more general functionals,
and we prove a uniqueness result for the
Monge-Kantorovich system (\ref{f:syst1})-(\ref{f:syst2}).
We also show some connection with the convergence
of solutions to the anisotropic $p$-Laplace equation
as $p\to\infty$.

\section{Notation and preliminaries}
\label{s:prel}

\subsection{Basic notation}
The standard scalar product of $x$ and $y\in\R^n$
will be denoted by $\pscal{x}{y}$, while $|x|$ will denote the
Euclidean norm of $x$.
Concerning the segment jointing $x$ with $y$, we set
\[
\cray{x,y} := \{tx+(1-t)y;\ t\in [0,1]\},
\qquad \oray{x,y} := \cray{x,y}\setminus\{x,y\}.
\]
As is customary, $B_r(x_0)$ and $\overline{B}_r(x_0)$
are respectively the open and the closed ball in $\R^n$
centered at $x_0$ and with radius $r>0$.

Given a set $A\subseteq \R^n$ its closure and its boundary
will be denoted by $\overline{A}$ and $\partial A$ respectively.

A bounded open set $\mathcal{O}\subset\R^n$
(or, equivalently,
$\overline{\mathcal{O}}$ or $\partial\mathcal{O}$)
is of class $C^k$, $k\in\N$,
if for every point $x_0\in\partial \mathcal{O}$
there exists a ball $B=B_r(x_0)$ and a one-to-one
mapping $\psi\colon B\to D$ such that
$\psi\in C^k(B)$, $\psi^{-1}\in C^k(D)$,
$\psi(B\cap \mathcal{O})\subseteq\{x\in\R^n;\ x_n > 0\}$,
$\psi(B\cap\partial \mathcal{O})\subseteq\{x\in\R^n;\ x_n = 0\}$.

If $f\colon\overline{\mathcal{O}}\to\R$ is measurable,
we define $\spt(f)$ as the intersection of
all closed sets $C\subseteq\overline{\mathcal{O}}$
such that $f\neq 0$ a.e.\ in $C$.


\subsection{Convex geometry}\label{ss:conv}
By $\convopen$ we shall denote the class of all
nonempty, compact, convex subsets of $\R^n$
with the origin as an interior point.
The polar set of $K\in\convopen$
is defined by
\[
K^0 = \{p\in\R^n;\ \pscal{p}{x}\leq 1\ \forall x\in K\}\,.
\]
We recall that, if $K\in\convopen$, then
$K^0\in\convopen$ and $K^{00} = (K^0)^0 = K$
(see \cite[Thm.~1.6.1]{Sch}).

Given $K\in\convopen$,
its gauge function is
\[
\gau{K}(\xi) = \inf\{ t\geq 0;\ \xi\in t K\}\,,
\quad \xi\in\R^n\,.
\]
Let $0<c_1\leq c_2$ be such that
$\overline{B}_{c_2^{-1}}(0)\subseteq K \subseteq \overline{B}_{c_1^{-1}}(0)$.
Upon observing that $\xi/\gau{K}(\xi)\in K$ for every $\xi\neq 0$,
we get
\begin{equation}\label{f:brho}
c_1 {|\xi|}\leq\gau{K}(\xi)\leq c_2 {|\xi|}\,,\quad
\forall\xi\in\R^n.
\end{equation}

We say that $K\in\convopen$ is of class $C^2_+$
if $\partial K$ is of class $C^2$
and all the principal curvatures are strictly positive
functions on $\partial K$.
We recall that, if $K$ is of class $C^2_+$, then $K^0$ is also of
class $C^2_+$ (see \cite[p.~111]{Sch}).

{}From now on we shall always assume that
\begin{equation}\label{f:ipoK}
K\in\convopen\
\textrm{is of class}\ C^2_+\,.
\end{equation}
Since $K$ will be kept fixed, from now on we shall use
the notation $\gauge = \gauge_K$ and $\pgauge=\gauge_{K^0}$.

We collect here some known properties
of $\gauge$ and $\pgauge$
that will be used
in the sequel (see e.g.\ \cite{CMg}, Theorem 2.1).

\begin{Theorem}\label{t:sch}
Let $K$ satisfy $(\ref{f:ipoK})$. Then the following hold:
\par\noindent (i)
The functions $\gauge$ and $\pgauge$ are
convex, positively $1$-homoge\-neous
in $\R^n$, and
of class $C^2$ in $\R^n\setminus\{0\}$.
\par\noindent (ii)
For every $\xi,\eta\in\R^n$, we have
\begin{equation}\label{f:sublin}
\gauge(\xi+\eta)\leq \gauge(\xi) + \gauge(\eta),\quad
\pgauge(\xi+\eta)\leq \pgauge(\xi) + \pgauge(\eta),
\end{equation}
and equality holds if and only if
$\xi = \lambda\, \eta$ or $\eta = \lambda\, \xi$
for some $\lambda\geq 0$.
\par\noindent (iii)
For every $\xi\neq 0$,
$D\gauge(\xi)$ belongs to $\partial K^0$ (i.e. $\pgauge(D\gauge(\xi))=1$), while $D\pgauge(\xi)$
belongs to $\partial K$ (i.e. $\gauge(D\pgauge(\xi))=1$). More precisely, $D\gauge(\xi)$
is the unique point of $\partial K^0$ such that
\[
\pscal{D\gauge(\xi)}{\xi} = \gauge(\xi),\ \textrm{and\ }
\pscal{p}{\xi} < \gauge(\xi)\ \forall p\in K^0,\ p\neq D\gauge(\xi)\,.
\]
Symmetrically,
the gradient of $D\pgauge(\xi)$ is the
unique point of $\partial K$ such that
\[
\pscal{D\pgauge(\xi)}{\xi} = \pgauge(\xi),\ \textrm{and\ }
\pscal{x}{\xi} < \pgauge(\xi)\ \forall x\in K,\ x\neq D\pgauge(\xi)\,.
\]
\end{Theorem}

\subsection{The distance function}\label{ss:dist}

Throughout the paper, we shall assume that
\begin{equation}\label{f:Omega}
\Omega\subset\R^n\
\textrm{ is a nonempty, bounded, open, connected set
of class $C^2$}.
\end{equation}

\begin{Definition}\label{d:dist}
Let $K\subset \R^n$ be a convex set fulfilling (\ref{f:ipoK}).
The Minkowski distance from the boundary of $\Omega$
associated to the convex body $K$
is the function $\dist\colon\overline{\Omega}\to\R$
defined by
\[
\dist(x) = \min_{y\in\partial\Omega}\pgauge(x-y),
\qquad
x\in\overline{\Omega}.
\]
\end{Definition}

Let us define the following spaces:
\begin{equation}\label{f:lipr}
\begin{split}
\Lipr & := \left\{
u\in C(\overline{\Omega});\
u(x)-u(y)\leq\pgauge(x-y)\
\forall \cray{x,y}\subset\overline{\Omega}
\right\}
\\ & \phantom{:}=
\left\{
u\in \Lip(\overline{\Omega});\
Du\in K\ \textrm{a.e.\ in}\ \Omega
\right\},
\end{split}
\end{equation}
\begin{equation}\label{f:liprz}
\Liprz := \left\{
u\in\Lipr;\
u(y) = 0\
\forall y\in\partial\Omega
\right\}.
\end{equation}
It is well known that $\dist\in\Liprz$, and that
$\dist\geq u$ for every $u\in\Liprz$
(see e.g.\ \cite{BaCD,Li}).

\begin{Definition}
The \textsl{inradius} of $\Omega$
is defined by $\inradius := \max\{\dist(x);\ x\in\Omega\}$.
\end{Definition}

We shall denote by $\proj(x)$ the set of projections of $x$
on $\partial\Omega$, that is
\begin{equation}\label{f:Pi}
\proj(x) := \{y\in\partial\Omega;\ \dist(x) = \pgauge(x-y)\},
\qquad x\in\overline{\Omega}.
\end{equation}
We recall that $x\mapsto \proj(x)$
is a sequentially upper semicontinuos multifunction, i.e.
\begin{equation}\label{f:uscpi}
x_k\in\overline{\Omega},\
y_k\in\proj(x_k),\ k\in\N;\
x_k\to x,\ y_k\to y
\Longrightarrow y\in\proj(x)\,.
\end{equation}

In some situations it will be convenient to consider
an extension $\dist^s$ of $\dist$
to $\R^n$ by setting
\[
\dist^s(x) = -\min_{z\in\overline{\Omega}} \pgauge(z-x),
\qquad x\in\R^n\setminus\Omega.
\]
This extension is the Minkowski signed distance from $\partial\Omega$.
Under the assumption (\ref{f:Omega}),
we have that $\dist^s$ is of class $C^2$ in a tubular
neighborhood $\mathcal{U}$ of $\partial\Omega$
(see \cite[Thm.~4.16]{CMf}).
In this neighborhood we can define the
Cahn-Hoffman vector field
\[
\ch(x) := -D\gauge(D\dist^s(x)),
\qquad
x\in\mathcal{U}.
\]
For every $y\in\partial\Omega$, the restriction
of $D\ch$ to the tangent space $T_y$ to $\partial\Omega$
at $y$ is a linear application from $T_y$ to $T_y$ having
$n-1$ real eigenvalues
$\curvg_1(y),\ldots,\curvg_{n-1}(y)$,
called the principal $\gauge$-curvatures or
anisotropic principal curvatures
of $\partial\Omega$ at $y$
(see \cite[Def.~5.5]{CMf}).
The anisotropic mean curvature is defined by
\begin{equation}\label{f:amc}
\mean(y) :=
\frac{1}{n-1}\sum_{i=1}^{n-1}\curvg_i(y)
=\dive\ch(y)
\qquad
\forall y\in\partial\Omega.
\end{equation}
A relevant quantity for the subsequent subjects will be
\begin{equation}\label{f:hze}
H_0 := \min\{\mean(y);\ y\in\partial\Omega\}\,.
\end{equation}

It can be shown that for every $x_0\in \Omega$ and $y\in\proj(x_0)$ the function
$\dist$ is differentiable at every $x\in \oray{y,x_0}$, and
\begin{equation}\label{f:gradd}
D\dist(x)=\dfrac{\nor(y)}{\gauge(\nor(y))}=D\pgauge(x_0-y),
\qquad \forall x\in\oray{y,x_0},\ y\in\proj(x_0),
\end{equation}
where $\nor(y)$ is the
Euclidean inward normal unit vector of $\partial\Omega$ at $y$.
Moreover one has
\[
\cray{y,x_0}=\{y+tD\gauge(\nor(y)),\ t\in [0, \dist(x_0)]\}\, ,\ y\in\proj(x_0)
\]
(see \cite{CMf}, Proposition 4.4).

\begin{Definition}
The singular set $\Sigma\subset\Omega$ of
$\dist$ is the set of all points in $\Omega$ where $\dist$
is not differentiable.
\end{Definition}
It is known that $x\in \Omega\setminus \Sigma$ if and only
if $x$ has a unique projection. Moreover
$\overline{\Sigma}$ has Lebesgue measure zero, and,
since $\dist$ is of class $C^2$ near the boundary,
$\overline{\Sigma}\subset\Omega$ (see \cite{CMf}, Corollary 6.9 and Theorem 4.16).

Notice that, thanks to (\ref{f:gradd}) and the positive $0$-homogeneity of
$D\gauge$, we infer that
\[
D\gauge(D\dist(x)) = D\gauge(\nor(y))\qquad
\forall x\in\Omega\setminus\Sigma,\
y\in\proj(x).
\]

\begin{Definition}
The normal distance of a point $x\in\overline{\Omega}$ to the cut locus
is defined by
\begin{equation}\label{f:tau}
\tau(x) :=
\begin{cases}
\min\{t\geq 0;\
x + t D\gauge(D\dist(x))\in\overline{\Sigma}\},
&\textrm{if $x\in\overline{\Omega}\setminus\overline{\Sigma}$},\\
0,
&\textrm{if $x\in\overline{\Sigma}$}.
\end{cases}
\end{equation}
The cut point $\cut(x)$ of
$x\in\overline{\Omega}\setminus\overline{\Sigma}$
is defined by
$\cut(x) := x + \tau(x)D\gauge(D\dist(x))$.
\end{Definition}

We recall that $\tau$ is a continuous function in $\overline{\Omega}$.
Furthermore, there exists $\mu > 0$ such that
$\tau(y) \geq \mu$ for every $y\in\partial\Omega$
(see \cite{CMf}, Lemma~4.1 and Theorem~6.7), and we have
\begin{equation}\label{f:tugl}
\tau(y)=\sup\{t\geq 0;\ y\in\Pi(y+tD\gauge(\nor(y)))\},\qquad \forall y\in\partial \Omega
\end{equation}
(see \cite{CMf}, Corollary 6.8).

\section{Existence}\label{s:exis}

In what follows we shall assume that
\begin{equation}\label{f:F}
\textrm{$f\in L^{\infty}(\Omega)$ is a nonnegative function}.
\end{equation}
Let us define the function
\begin{equation}\label{f:vf}
\vf(x) :=
\begin{cases}
\displaystyle{
\int_0^{\tau(x)} f(\Phi(x,t))
\prod_{i=1}^{n-1}\frac{1-(\dist(x)+t)\, \curvg_i(x)}%
{1-\dist(x)\, \curvg_i(x)}\, dt}
&\textrm{if $x\in{\overline{\Omega}}\setminus\overline{\Sigma}$},\\
0,
&\textrm{if $x\in\overline{\Sigma}$}\,,
\end{cases}
\end{equation}
where, for $x\in\overline{\Omega}\setminus\overline{\Sigma}$
and $\proj(x) = \{y\}$, we have set
\[
\Phi(x, t) := x + t\, D\gauge(D\dist(x)),
\qquad
\curvg_i(x) := \curvg_i(y),\quad
i=1,\ldots,n-1\,.
\]
Since the maps $\tau$ and $\curvg_i$, $i=1,\ldots,n-1$,
are continuous in $\overline{\Omega}\setminus\overline{\Sigma}$,
and the map $\Phi$ is continuous in
$(\overline{\Omega}\setminus\overline{\Sigma})\times\R$
(see \cite[Section~7]{CMf}), the function $\vf$ is well defined
and bounded in $\Omega$.

Let $c\colon\R\times (0,\infty)\to\R$ be the function defined by
\begin{equation}\label{f:ipr}
c(t,r) :=
\begin{cases}
\frac{1-(1-t\, r)^n}{n t}\,,
&\textrm{if $t\neq 0$},\\
r\,,
&\textrm{if $t = 0$}.\\
\end{cases}
\end{equation}
It is straightforward to check that
$t\mapsto c(t,r)$ is a strictly monotone decreasing function for
$t\leq \frac{1}{r}$.

The main properties of the function $\vf$ defined in (\ref{f:vf}) are the following.
\begin{Proposition}\label{p:vf}
The function
$\vf$ belongs to $L^{\infty}(\Omega)$, $\vf\geq 0$ a.e.\ in $\Omega$, and
\begin{equation}\label{f:vfestii}
{\|\vf\|}_{L^\infty(Q)}< {\|f\|}_{\infty}\,c(H_0, \inradius),
\end{equation}
for every compact set $Q\subset\Omega$,
where $H_0$ is the constant defined in (\ref{f:hze}).
In particular
\begin{equation}\label{f:vfesti}
{\|\vf\|}_{\infty}\leq {\|f\|}_{\infty}\,c(H_0, \inradius).
\end{equation}
Moreover $\vf$ satisfies
\begin{equation}\label{f:vfeq}
\int_{\Omega} \vf\pscal{D\gauge(D\dist)}{D\varphi}\, dx
= \int_{\Omega}f\, \varphi\, dx
\qquad
\forall\varphi\in W^{1,1}_0(\Omega).
\end{equation}
\end{Proposition}

\begin{proof}
The facts that $\vf\in L^{\infty}(\Omega)$
and $\vf$ is a solution to
(\ref{f:vfeq})
are proved in \cite[Section~7]{CMf}
in the case $f\in C(\overline{\Omega})$.
The very same proof also works for $f\in L^{\infty}(\Omega)$.
The estimates (\ref{f:vfestii}) and (\ref{f:vfesti}) follows from a
straightforward adaptation of Proposition~5.9 in \cite{CCCG}.
\end{proof}

\begin{Theorem}\label{t:exif}
Let $\Omega\subset\R^n$ be a set fulfilling  (\ref{f:Omega}), and
let $f$ satisfy (\ref{f:F}).
Assume that (H1)--(H2) hold, together with
\begin{itemize}
\item[(H3)]
there exists
$H^*\leq H_0$
such that
$c(H^*,\inradius)\, {\|f\|}_{\infty}\leq \Lambda,$
\end{itemize}
where $\Lambda\in [0,+\infty]$ is the constant defined
in (\ref{lak}).
Then the distance function $\dist$
is a minimizer of $J$ in $W^{1,1}_0(\Omega)$.
Moreover, if $u$ is a
minimizer of $J$ in $W^{1,1}_0(\Omega)$, then
$u\in\Liprz$.
\end{Theorem}

\begin{proof}
Since $H_0\leq 1/\inradius$ (see \cite{CMf}, Lemma 5.4), we have that
$c(H^*,\inradius)\geq c(H_0,\inradius)$.
{}From (H3) and (\ref{f:vfesti}) we infer that
${\|\vf\|}_{\infty}\leq \Lambda$.
Hence, if $u\in W^{1,1}_0(\Omega)$ we have that
\begin{equation}\label{f:estih}
\begin{split}
J(u)-J(\dist) & \geq
\int_\Omega\left[\max\{\Lambda(\gauge(Du)-1),\, 0\}-f(u-\dist)\right]\, dx
\\ & \geq\int_\Omega\left[\vf(\gauge(Du)-1)-f(u-\dist)\right]\, dx
\\ & \geq \int_\Omega\left[\vf\pscal{D\gauge(D\dist)}{Du-D\dist}-f(u-\dist)\right]\, dx=0\,,
\end{split}
\end{equation}
where the third inequality is
a consequence of
\[
\pscal{D\gauge(p)}{q-p} =
\pscal{D\gauge(p)}{q} - 1
\leq \gauge(q) -1
\quad
\forall p\in \partial K,\ q\in\R^n
\]
(see Theorem \ref{t:sch}(iii)), while the last equality follows from (\ref{f:vfeq}).

Moreover, by (\ref{f:vfestii})
the strict inequality holds in (\ref{f:estih})
whenever $Du(x)\not\in K$ on a set of positive Lebesgue measure.
Hence $J(u) = J(\dist)$ implies
$Du\in K$ a.e.\ in $\Omega$.
\end{proof}

\begin{Remark}
A more precise existence result,
containing detailed information about the structure of
minimizers,
will be given in Theorem~\ref{t:uf} below.
\end{Remark}

\begin{Remark}
If $\Omega$ is a convex set, then $H_0\geq 0$,
and $c(H_0,\inradius)\leq\inradius$.
Hence (H3)
is certainly satisfied if $\inradius {\|f\|}_{\infty}\leq \Lambda$.
This point
was already highlighted by A.~Cellina in \cite{Ce1}.
\end{Remark}

\begin{Remark}
The assumption (H3) can be read as a growth condition on $h$
outside $K$, which guarantees that $J$ is bounded from below.
In \cite{CCCG}, Example~5.6, it is exhibited a functional $J$,
not bounded from below,
which satisfies all the assumptions of Theorem~\ref{t:exif} but (H3).
\end{Remark}

\section{Uniqueness}
\label{s:uni}

The aim of this section is to prove the
following uniqueness result.

\begin{Theorem}\label{t:unif}
Let $\Omega\subset\R^n$ be a set fulfilling  (\ref{f:Omega}).
Assume that (H1)--(H3) hold.
Then $\dist$ is the unique minimizer of $J$
in $W^{1,1}_0(\Omega)$
if and only if
$\Sigma\subset\spt(f)$.
\end{Theorem}

\begin{Remark}
Although hypothesis (\ref{f:Omega}) is needed for proving the uniqueness result,
the preliminary
results in Lemma~\ref{l:cran},
Proposition~\ref{l:us}, Lemma~\ref{l:ray}, and Proposition~\ref{p:closed}, having an interest by themselves,
are proved without using the regularity
assumption on $\Omega$.
\end{Remark}

The following result is essentially due to
M.G.\ Crandall \cite[Lemma~7.3]{Cran}
(see also \cite[Prop.~4.2]{Amb}). We have to make
some minor changes with respect to the Euclidean case,
due to the fact that the function $\pgauge$ need not be
symmetric.

\begin{Lemma}\label{l:cran}
Let  $u\in\Lipr$, $\oray{x_0, x_1}\subset {\Omega}$ ($x_0\neq x_1$),
and assume that
\begin{equation}\label{f:ulin}
u(x_0+t(x_1-x_0)) = u(x_0) + t\,\pgauge(x_1-x_0)
\qquad\forall t\in [0,1].
\end{equation}
Then $u$ is differentiable at every point $x\in \oray{x_0,x_1}$,
and $Du(x) = D\pgauge(p)$, where
$p := (x_1-x_0)/\pgauge(x_1-x_0)$.
\end{Lemma}

\begin{proof}
Let $\overline{x}\in \oray{x_0,x_1}$.
By (\ref{f:ulin}) there exists $\varepsilon>0$ such that $\overline{x}+\sigma p \in \Omega$
and $u(\overline{x}+\sigma\, p)=u(\overline{x})+\sigma$ for every $\sigma\in(-\varepsilon, \varepsilon)$.
It is not restrictive to assume that
$\overline{x}=0$ and $u(\overline{x}) = 0$. Hence we have
\begin{equation}\label{f:ulinx}
\sigma\, p \in \Omega,\ u(\sigma\, p) = \sigma,\qquad \forall \sigma\in (-\varepsilon, \varepsilon)\,.
\end{equation}
For every $x\in\R^n$
let us define $Px := \pscal{D\pgauge(p)}{x}p$, and let $0<\delta<\varepsilon$ be such that
\[
|\pscal{D\pgauge(p)}{x}+s|<\varepsilon,\qquad \forall x\in B_{\delta}(0)\subseteq \Omega,\
\forall s\in (-\delta,\delta).
\]
Fixed $r\in (0,\delta)$,
for every $x\in B_{\delta}(0)$, the point $y=Px+rp$ belongs to
$\Omega$. Moreover, from (\ref{f:ulinx}) and (\ref{f:lipr})
we obtain
\[
\begin{split}
\pscal{D\pgauge(p)}{x}+r & -\pgauge\left(Px+r\,p-x\right)
\\ & =
u\left(Px+r\,p\right)-\pgauge\left(Px+r\,p-x\right)
\leq u(x)
\end{split}
\]
so that
\begin{equation}\label{f:ulin1}
-r\left[\pgauge\left(\frac{Px-x}{r}+p\right)-1\right]
\leq u(x)- \pscal{D\pgauge(p)}{x}\,.
\end{equation}
Similarly, the point $z=Px-rp$ belongs to $\Omega$, and
\begin{equation}\label{f:ulin2}
u(x)- \pscal{D\pgauge(p)}{x}\leq
r\left[\pgauge\left(\frac{x-Px}{r}+p\right)-1\right]\,.
\end{equation}
{}From the differentiability of $\pgauge$ at $p$, and the fact that
$\pgauge(p) = 1$,
we have
\[
\pgauge(p+q) = \pgauge(p)+\pscal{D\pgauge(p)}{q}+o(|q|)
= 1 + \pscal{D\pgauge(p)}{q}+o(|q|), \quad q\to 0\,.
\]
Combining (\ref{f:ulin1}) and (\ref{f:ulin2}), and using the fact that
$\pscal{D\pgauge(p)}{x-Px} = 0$, we get
$u(x)- \pscal{D\pgauge(p)}{x}= o(|x|)$, $x\to 0$.
Hence $u$ is differentiable at $\overline{x}=0$
and $Du(\overline{x}) = D\pgauge(p)$.
\end{proof}
In Section~\ref{s:exis} we have proved that $\dist$ is a minimizer of $J$ in
$W^{1,1}_0(\Omega)$.
Since every other minimizer has to belong to $\Liprz$, $\dist$
turns out to be the maximal minimizer of $J$.
Now we want to construct the minimal
minimizer of $J$ (see (\ref{f:uf}) below).
This minimizer will be obtained
as the supremum of cone--shaped functions.

Assume that
\begin{equation}\label{f:S}
\textrm{$S$ is a nonempty closed subset of $\overline{\Omega}$},
\end{equation}
and define the function
\begin{equation}\label{f:maxz}
u_S(x) := \max_{z\in S\cup\partial\Omega} [\dist(z)-\pgauge(z-x)],
\qquad x\in\overline{\Omega}\,.
\end{equation}

The \textsl{transport set} of $S$ is defined by
\begin{equation}\label{f:ts}
\TS :=
\bigcup_{\substack{z\in S\\ y\in\proj(z)}}
 \oray{y,z}
\subset{\Omega}\,.
\end{equation}
The main properties of the function $u_S$ are collected in the following result.

\begin{Proposition}\label{l:us}
Let $S\subseteq\overline{\Omega}$ satisfy (\ref{f:S}) and
let $u_S$ be the function defined in (\ref{f:maxz}).
Then the following hold.
\begin{itemize}
\item[(i)]
$u_S\in\Liprz$.
\item[(ii)]
$u_S \leq \dist$ in $\overline{\Omega}$, and
$u_S = \dist$ in $\overline{\TS}\cup S$.
Moreover, for every $z\in {S}\cap\Omega$ and
$y\in\proj(z)$,
$u_S$ is differentiable at every $x\in\oray{y,z}$
and $Du_S(x) = D\dist(x)$.
\item[(iii)]
If $u\in\Liprz$ and $u=\dist$ in $S$,
then $u\geq u_S$ in $\Omega$.
\item[(iv)]
$\gauge(Du_S(x)) = 1$
at every $x\in\Omega\setminus {S}$ where
$u_S$ is differentiable.
\item[(v)]
$u_S = \dist$ in $\Omega$ if and only if
$\Sigma\subseteq {S}$.
\end{itemize}
\end{Proposition}

\begin{proof}
(i)
It is plain that $u_S$ is a Lipschitz function vanishing on $\partial\Omega$.
It remains to prove that
\begin{equation}\label{f:us1}
u_S(x_1)-u_S(x_2)\leq\pgauge(x_1-x_2)
\qquad\forall x_1,x_2\in {\Omega}.
\end{equation}
Let $x_1,x_2\in {\Omega}$,
and let $z\in S\cup\partial\Omega$ be such that
\begin{equation}\label{f:us2}
u_S(x_1)=\dist(z)-\pgauge(z-x_1).
\end{equation}
By the very definition of $u_S$ we obtain
\[
-u_S(x_2)\leq -\dist(z)+\pgauge(z-x_2)
\leq -\dist(z) + \pgauge(z-x_1) + \pgauge(x_1-x_2),
\]
which, combined with (\ref{f:us2}), gives (\ref{f:us1}).

\smallskip
\noindent
(ii)
The inequality $u_S\leq\dist$ follows from the fact that
$u_S\in\Liprz$.
If $z\in S$, we have that
$u_S(z) \geq \dist(z)$, hence
$u_S(z) = \dist(z)$.

In order to prove that $u_S = \dist$ in $\overline{\TS}$,
let us fix $z\in S\cap\Omega$, $y\in\proj(z)$,
and let $x\in \oray{y,z}$, so that $\proj(x)=\{y\}$.
{}From Theorem \ref{t:sch}(iii) we infer that
\[
\begin{split}
\pgauge(x-y) & = \dist(x)\geq u_S(x)\geq\dist(z)-\pgauge(z-x)
\\ & =
\pgauge(z-y)-\pgauge(z-x)
= \pgauge(x-y).
\end{split}
\]
Hence  $\dist = u_S$ in $\TS$, and the equality extends
by continuity to $\overline{\TS}$.

Finally,
we have that
\[
u_S(y+t(z-y)) = \dist(y+t(z-y)) = t\, \pgauge(z-y),
\qquad\forall t\in [0,1].
\]
{}From Lemma~\ref{l:cran}
we conclude that
$u_S$ is differentiable at every $x\in\oray{y,z}$,
and $D u_S(x) = D\pgauge(z-y) = D\dist(x)$
(see (\ref{f:gradd})).

\smallskip
\noindent
(iii)
Observe that $u=\dist$ on $S\cup\partial\Omega$.
Let $x_0\in\Omega\setminus {S}$,
and let us prove that $u(x_0)\geq u_S(x_0)$.
Let $z\in S\cup\partial\Omega$ be such that
$u_S(x_0) = \dist(z)-\pgauge(z-x_0)$.
Then we have that
\[
u_S(x_0)=\dist(z)-\pgauge(z-x_0) = u(z)-\pgauge(z-x_0) \leq u(x_0)\,,
\]
where the last inequality follows from the fact that $u\in \Lipr$.

\smallskip
\noindent
(iv)
{}From the general theory of marginal functions we have that for every $x\in \Omega\setminus S$
where $u_S$ is differentiable there exists $z\in S\cup\partial\Omega$ such that
$Du_S(x)=D\pgauge(z-x)$ (see e.g.\ \cite{CaSi}, Theorem~3.4.4).
Then the conclusion follows from Theorem \ref{t:sch}(iii).

\smallskip
\noindent
(v)
Assume that $\Sigma\subseteq {S}$.
{}From (ii) we have that
$u_S(x) = \dist(x)$ for every $x\in\overline{\Sigma}$.
Let $x\in\Omega\setminus\overline{\Sigma}$.
Let $y$ be the unique projection of $x$ on $\partial\Omega$
and let $z\in \overline{\Sigma}$
be the cut point of $x$.
Then we have
\[
\dist(x)=\pgauge(x-y)
= \pgauge(z-y)-\pgauge(z-x)
\leq u_S(x)\leq\dist(x),
\]
hence we conclude that $u_S(x) = \dist(x)$.

Conversely, assume that $u_S = \dist$ in $\overline{\Omega}$.
By contradiction, suppose that there exists
a point $x_0\in\Sigma$, $x_0\not\in S$.
Let $z\in S\cup\partial\Omega$ be such that
$u_S(x_0) = \dist(z) - \pgauge(z-x_0)$.
Since $u_S(x_0)=\dist(x_0)>0$,
it is plain that $z\in S\cap\Omega$, for otherwise we would have
$u_S(x_0)\leq 0$.
Summarizing, we have that
\begin{equation}\label{f:lus0}
\exists z\in S\cap\Omega:\quad\dist(x_0) = u_S(x_0)
= \dist(z) - \pgauge(z-x_0).
\end{equation}
On the other hand, we are going to show that
the assumptions $x_0\in\Sigma$, $x_0\neq z$, imply that
\begin{equation}\label{f:lus1}
\dist(z) < \dist(x_0) + \pgauge(z-x_0).
\end{equation}
Namely, let $y_0\in\proj(x_0)$.
If $y_0\in\proj(z)$, then
$x_0\not\in \oray{y_0, z}$, since otherwise $\dist$ should be
differentiable at $x_0$, in contrast with the
assumption $x_0\in\Sigma$.
Moreover, from (\ref{f:lus0}), $\dist(x_0)<\dist(z)$,
so that $z\not\in\oray{y_0,x_0}$.
Hence the three points $y_0$, $x_0$, $z$ do not lie on the
same ray, so that
\[
\dist(z) = \pgauge(z-y_0) < \pgauge(z-x_0)+\pgauge(x_0-y_0)
= \pgauge(z-x_0)+\dist(x_0),
\]
and (\ref{f:lus1}) holds.
On the other hand, if
$y_0\not\in\proj(z)$, then
\[
\dist(z) < \pgauge(z-y_0) \leq \pgauge(z-x_0)+\pgauge(z_0-y_0)
= \pgauge(z-x_0)+\dist(x_0),
\]
and again (\ref{f:lus1}) holds.
\end{proof}


\begin{Definition}\label{d:ss}
The \textsl{reduced set} $S^*$
is the set of points
$z^*\in S$ such that the following holds:
if $z\in S$, $y\in\proj(z)$ and $z^*\in\cray{y,z}$,
then $z=z^*$.
\end{Definition}

\begin{Lemma}\label{l:ray}
Let $z\in S$ and $y\in\proj(z)$.
Then there exists a unique $z^*\in S^*$ such that
$\cray{y,z}\subseteq\cray{y,z^*}$.
\end{Lemma}

\begin{proof}
The uniqueness of $z^*$ is a straightforward
consequence of the definition of $S^*$.

For what concerns the existence,
assume first that $z\in\Omega$, so that
$y\neq z$.
Let us define
\[
A :=
\left\{
\lambda\geq 0;\
y_{\lambda} := y + \lambda\, \frac{z-y}{\pgauge(z-y)}\in S,\
y\in\proj(y_{\lambda})
\right\}\,.
\]
We have that $A$ is bounded, 
$A$ is closed
(see (\ref{f:uscpi})),
and $\dist(z)\in A$, hence $A$
admits a maximum $\sigma\geq\dist(z)$.

We claim that the point $z^* := y_{\sigma}$
belongs to $S^*$.
Namely, given $z'\in S$, $y'\in\proj(z')$
such that $z^*\in\cray{y',z'}$,
we have to prove that $z'=z^*$.
If $y'=y$, then the equality follows from the maximality
of $\sigma$.
On the other hand, if $y'\neq y$, we have that
$\proj(z^*)$ contains two different points $y$ and $y'$,
so that $z^*\in\Sigma$.
Since $\dist$ is differentiable at every point in $\oray{y',z'}$, we infer that $z^*=z'$.

Consider now the case $z\in\partial\Omega\cap S$,
so that $y=z$.
If $z\not\in\proj(z')$ for every $z'\in\Omega\cap S$,
then $z\in S^*$.
Otherwise let $z'\in\Omega\cap S$ be such that
$y=z\in\proj(z')$.
{}From the first part of the proof,
there exists $z^*\in S^*$ such that
$\{z\} = \cray{y,z}\subset\cray{y,z'}\subseteq\cray{y,z^*}$,
and the proof is complete.
\end{proof}

Each segment $\oray{y,z}$, with $z\in S^*$ and $y\in\proj(z)$,
will be called a \textsl{transport ray}.
We shall denote by $\TSe$ the union of the closures
$\cray{y,z}$ of all transport rays, that is
\begin{equation}\label{f:tse}
\TSe :=
\bigcup_{\substack{z\in S^*\\ y\in\proj(z)}}
 \cray{y,z}
\subset\overline{\Omega}
\end{equation}
(To be precise, $\TSe$ also contains the points
$z\in S^*\cap\partial\Omega$,
that are not closures of transport rays.)

{}From the definition it is plain that
two different transport rays have empty intersection, and that
transport rays are maximal, i.e.\ if
$\oray{y,z}$ is a transport ray,
then $\oray{y,y+t(z-y)}$ is not a transport ray
for every $t>1$.

\begin{Proposition}\label{p:closed}
Let $S$ satisfy (\ref{f:S}).
Then $\TSe$ is a closed subset of $\overline{\Omega}$
and
\begin{equation}\label{f:tt}
\TSe = \overline{\TS}\cup S =
\bigcup_{\substack{z\in S\\ y\in\proj(z)}}
 \cray{y,z}\,.
\end{equation}
\end{Proposition}

\begin{proof}
Let $T$ denote the set appearing in the
right-hand side of (\ref{f:tt}).
One easily obtain that $T=\TSe$. Namely,
since $S^*\subseteq S$,
we have that $\TSe\subseteq T$, while
the inclusion $T\subseteq\TSe$ follows from
Lemma~\ref{l:ray}.

It remains to prove that $\TSe=\overline{\TS}\cup S$.
Let us start by proving that $\TSe$ is a closed subset of
$\overline{\Omega}$.
Let $(x_j)\subset\TSe$ be a sequence
converging to a point $x$.
By definition of $\TSe$,
for every $j\in\N$ there exist $z_j\in S^*$, $y_j\in\proj(z_j)$
and $t_j\in [0,1]$ such that
$x_j = y_j+t_j(z_j-y_j)$.
We can extract a subsequence (not relabeled) such that
$z_j\to z\in S$, $y_j\to y\in\partial\Omega$,
$t_j\to t\in [0,1]$,
so that $x=y+t(z-y)$.
By the upper semicontinuity of the multifunction $\proj$,
we have that $y\in\proj(z)$,
so that
$x = y+t(z-y)\in\cray{y,z}$
with $z\in S$ and $y\in\proj(z)$, i.e.
$x\in T = \TSe$.

Finally, in order to prove that $\TSe = \overline{\TS}\cup S$
it is enough to observe that
$\TS\cup S\subseteq\TSe$ and $\TSe = T \subseteq \overline{\TS}\cup S$.
\end{proof}


Let us define the function
\begin{equation}\label{f:lstar}
\lray(x) :=
\begin{cases}
\pgauge(z-y),
&\textrm{if $x\in\cray{y,z}$ for some $z\in S^*$ and $y\in\proj(z)$},\\
0,
&\textrm{otherwise in $\overline{\Omega}$}.
\end{cases}
\end{equation}

\begin{Proposition}[Upper semicontinuity of $\lray$]\label{p:usc}
Let $S$ satisfy (\ref{f:S}). Then
the function $\lray$ is upper semicontinuous in $\overline{\Omega}$.
\end{Proposition}

\begin{proof}
Let $(x_j)\subset\overline{\Omega}$ be a sequence
converging to a point $x\in\overline{\Omega}$.
We have to prove that
$\lray(x)\geq\limsup_j\lray(x_j)$.
It is not restrictive to assume that
$\limsup_j\lray(x_j)=\lim_j\lray(x_j)$.
If $\lim_j\lray(x_j) = 0$ the conclusion is trivial,
hence it is enough to consider only the case
$(x_j)\subset\TSe$.
By Proposition~\ref{p:closed}, we have that also $x\in\TSe$.
By definition, for every $j\in\N$ there exist
$z_j\in S^*$, $y_j\in\proj(z_j)$ and
$t_j\in [0,1]$ such that
$x_j = y_j+t_j\, (z_j-y_j)$,
so that $\lray(x_j) = \pgauge(z_j-y_j)$.
We can pass to a subsequence (not relabeled)
such that
$z_j\to z\in S$, $y_j\to y\in\proj(z)$,
$t_j\to t\in [0,1]$,
so that $x = y+t(z-y)$.
{}From Lemma~\ref{l:ray} there exists $z^*\in S^*$
such that $\cray{y,z}\subseteq\cray{y,z^*}$,
hence
\[
\lray(x) = \pgauge(z^*-y)\geq
\pgauge(z-y)=
\lim_j\pgauge(z_j-y_j)=\lim_j\lray(x_j),
\]
concluding the proof.
\end{proof}

The following proposition says that
the measure of the set $S^*$ of the endpoints of
transport rays has zero Lebesgue measure.
This is a well known property in transport theory
(see e.g.\ \cite[Coroll.~6.1]{Amb}
or \cite[Lemma~2.15]{FMC}).
We give a direct proof based on the techniques
developed in \cite{CMf}
(see in particular Corollary~4.15 in \cite{CMf}).
Here $\meas$ denotes the $n$-dimensional Lebesgue measure.

\begin{Proposition}\label{p:vanish}
The reduced set $S^*$ has vanishing Lebesgue measure.
As a consequence,
$\meas(\TSe\setminus\TS) = 0$.
\end{Proposition}

\begin{proof}
By the very definition of $S^*$ we have that
\[
S^*\subset Z :=
\{y+\lray(y)\, D\gauge(\nor(y));\
y\in\partial\Omega
\}.
\]
Let $Y_k\colon\U_k\to\R^n$, $\U_k\subset\R^{n-1}$ open,
$k=1,\ldots,N$,
be local parameterizations of $\partial\Omega$ of class $C^2$,
such that $\bigcup_{k=1}^{N} Y_k(\U_k) = \partial\Omega$.
For every $k=1,\ldots,N$,
let $\Psi_k\colon \U_k\times\R\to\R^n$ be the map
\[
\Psi_k(y,t) = Y_k(y) + t D\gauge(\nor(Y_k(y))),
\quad (y,t)\in \U_k\times\R\,.
\]
For every $k=1,\ldots,n$ let $U_k\subset\U_k$
be a compact set such that $\bigcup_k Y_k(U_k)$ covers
$\partial\Omega$,
and let
\[
A_k = \{(y,t);\ y\in U_k,\ t\in [0,\lray(Y_k(y))]\}\,.
\]
{}From Proposition~\ref{p:usc},
$\lray$ is an upper semicontinuous function,
hence for every $k=1,\ldots,n$,
$A_k$ is a compact set and
the Lebesgue measure of
the graph
\[
\Psi_k^{-1}(Z)\cap A_k =
\{(y,t)\in U_k\times\R;\ t=\lray(Y_k(y))\}
\]
vanishes.
Moreover $\Psi_k\in C^1(\U_k\times\R)$ for every $k=1,\ldots,N$
(see \cite{CMf}, Theorem~4.13),
hence $\Psi_k$ is Lipschitz continuous on
the compact set $A_k$.
Let $L$ be the maximum of the Lipschitz constants
of the functions $\Psi_1,\ldots,\Psi_N$.
Since $\bigcup_{k=1}^N \Psi_k(A_k) = \overline{\Omega}$,
and hence
$S^*\subseteq Z\subseteq\bigcup_{k=1}^N \Psi_k\left(\Psi_k^{-1}(Z)\cap A_k\right)$,
we finally get
\[
\meas(S^*)\leq
\sum_{k=1}^N \meas\left[\Psi_k\left(
\Psi_k^{-1}(Z)\cap A_k\right)\right]
\leq L^n
\sum_{k=1}^N \meas\left[
\Psi_k^{-1}(Z)\cap A_k\right]
= 0\,.
\]
Finally, the last assertion follows from
the inclusion
$\TSe\setminus\TS\subseteq\partial\Omega\cup S^*$
(see Lemma~\ref{l:ray}).
\end{proof}

\begin{Remark}
In general $S^*$ is not closed.
For example, let $\Omega = B_R(0)\subset\R^2$, with $R>2\pi+1$, and let
\[
S := \{(r\, \cos\theta, r\, \sin\theta);\
\theta\in [0,2\pi),\ \theta+1\leq r\leq R\}.
\]
It is easy to verify that $S$ is a closed set.
On the other hand, if we consider the Euclidean metric ($K=\overline{B}_1(0)$),
we have that the reduced set is
\[
S^* = \{((\theta+1)\,\cos\theta, (\theta+1)\,\sin\theta);\
\theta\in [0,2\pi)\},
\]
which is not closed since $\overline{S^*}\setminus S^* = \{(2\pi+1,0)\}$.
We can also easily compute the function $\lambda^*$,
\[
\lambda^*(r\cos\theta, r\sin\theta) = R-\theta-1,\quad
\theta+1\leq r \leq R,\ 0\leq\theta < 2\pi,
\qquad \lambda^* = 0\ \textrm{otherwise in}\ \Omega,
\]
which is upper semicontinuous but not continuous in $\overline{\Omega}$.
\end{Remark}

For future reference, we collect here the properties
of the function $u_S$ defined in (\ref{f:maxz}).

\begin{Corollary}\label{c:us}
Let $\Omega$ and $S$ satisfy (\ref{f:Omega})
and (\ref{f:S}) respectively.
Then the function $u_S$
belongs to $\Liprz$ and satisfies the following
properties.
\begin{itemize}
\item[(i)]
$u_S = \dist$ on the closed set
$\TSe\supseteq S$, and
$D u_S(x) = D\dist(x)$ for a.e.\ $x\in\TSe$.
\item[(ii)]
$\gauge(Du_S) = 1$ a.e.\ in $\Omega$.
\item[(iii)] If $u\in\Liprz$ coincides with
$\dist$ on $S$, then $u_S\leq u\leq\dist$ in $\overline{\Omega}$,
and $Du(x) = D\dist(x)$ for a.e.\ $x\in\TSe$.
\item[(iv)]
$u_S = \dist$ if and only if $\Sigma\subseteq S$.
\end{itemize}
\end{Corollary}

\begin{proof}
(i) From Proposition~\ref{p:closed} we have that
$\TSe = \overline{\TS}\cup S$.
Hence, from Proposition~\ref{l:us}(ii), we deduce that
$u_S = \dist$ on $\TSe$ and $u_S$ is differentiable
at every point of $\TS$, with $D u_S = D\dist$.
Since, by Proposition~\ref{p:vanish}, $\meas(\TSe\setminus\TS) = 0$,
we conclude that $u_S$ is differentiable a.e.\ on $\TSe$ with
$D u_S = D\dist$.

(ii)
{}From Proposition~\ref{l:us}(iv) we have
that $u_S$ is differentiable and $\gauge(D u_S) = 1$
a.e.\ on $\Omega\setminus S$.
On the other hand, by (i) we have that $u_S$ is differentiable
for a.e.\ $x\in\TSe\supseteq S$, and
$\gauge(D u_S(x)) = \gauge(D\dist(x)) = 1$.

(iii) and (iv) are proved in
Proposition~\ref{l:us} (iii) and (v) respectively,
with the exception of the equality $Du = D\dist$ a.e.\ on $\TSe$,
which follows from Lemma~\ref{l:cran} upon observing that
$u$ coincides with $\dist$ along transport rays.
\end{proof}

\begin{Theorem}\label{t:uf}
Let the assumptions of Theorem~\ref{t:exif} hold.
Then the function
\begin{equation}\label{f:uf}
u_f(x) := \max_{z\in\spt(f)\cup\partial\Omega}[\dist(z)-\pgauge(z-x)]\,,
\qquad x\in\overline{\Omega}
\end{equation}
is a minimizer of $J$ in $W^{1,1}_0(\Omega)$.
Moreover, any other minimizer $u$ of $J$ in $W^{1,1}_0(\Omega)$
belongs to $\Liprz$ and satisfies $u_f\leq u\leq d$.
In particular, $u=\dist$ on the set
$\Tf := T^e_{\spt(f)}$ defined in (\ref{f:tse}) with $S=\spt(f)$.
\end{Theorem}

\begin{proof}
Observe that $u_f$ is the function defined
by (\ref{f:maxz}) with $S = \spt(f)$.
Let $\vf$ be the function defined in (\ref{f:vf}).

We claim that
\begin{equation}\label{f:vanv}
\vf(x)=0\,, \qquad \forall x\in\Omega\setminus\Tf\,.
\end{equation}
Property (\ref{f:vanv}) is plain if $x\in\overline{\Sigma}\setminus\Tf$.
Assume now that
$x\in \Omega\setminus\Tf$, $x\not\in \overline{\Sigma}$, and let $z\in\overline{\Sigma}$
be its cut point. Since $\spt(f)\subseteq\Tf$, in order to show that $\vf(x) = 0$
it is enough to prove that $\oray{x,z}$ does not intersect $\Tf$.
Assume, by contradiction,
that there exists $x_0\in \oray{x,z}\cap \Tf$.
Then $x_0\not \in \Sigma$, and $\Pi(x_0)=\Pi(x)=:\{y\}$.
Moreover, since $x_0\in \Tf$, then $\cray{y,x_0}\subseteq \Tf$,
which contradicts the fact that $x\not\in \Tf$.

By Corollary~\ref{c:us}(i) and (\ref{f:vanv}),
we obtain that
$u_f=\dist$  and
$D u_f=D\dist$ a.e.\ in the set $\{\vf > 0\}$.
Hence,
by (\ref{lak}), (\ref{f:vfesti}), (H3),
and Corollary~\ref{c:us}(ii), we conclude that
\[
J(u) - J(u_f)\geq
\int_{\Omega}
\left[\vf\pscal{D\gauge(D\dist)}{Du-D u_f}-(u-u_f)f\right]\, dx
\]
for every $u\in W^{1,1}_0(\Omega)$
(see (\ref{f:estih})).
Then, by (\ref{f:vfeq}), we get
$J(u)\geq J(u_f)$,
i.e.\ $u_f$ is a minimizer of $J$.

Assume now that
$u$ is a minimizer of $J$ in $W^{1,1}_0(\Omega)$.
Since $\dist$ is also a minimizer,
and $u\in\Liprz$ (see Theorem~\ref{t:exif}),
we have that
\[
0=J(u)-J(\dist) = \int_{\Omega}\left[h(Du)+(\dist-u)f\right]\, dx.
\]
Since $h$, $f$ and $\dist-u$ are non-negative functions,
it follows that $h(Du) = 0$ a.e.\ in $\Omega$,
and $u=\dist$ in $\spt(f)$.
{}From Corollary~\ref{c:us}(iii)
it follows that $u_f\leq u\leq d$ in $\Omega$ and $u=\dist$ on $\Tf$.
\end{proof}

\begin{Remark}
It is not difficult to see that
the conclusions of Theorem~\ref{t:uf} continue to hold
under the following assumptions:
$\Omega$ and $K$ satisfy (\ref{f:Omega}) and
(\ref{f:ipoK}) respectively,
$f\in L^1(\Omega)$ is a non-negative function,
and $h$ is the indicator function of the set $K$.
\end{Remark}

\begin{proof}[Proof of Theorem~\ref{t:unif}]
{}From Theorem~\ref{t:uf}
we have that $u_f$ and $\dist$ are both minimizers of
$J$ in $W^{1,1}_0(\Omega)$,
and that any other minimizer $u$ of $J$ satisfies
$u_f\leq u\leq\dist$.
Hence $\dist$ is the unique minimizer of $J$
if and only if $u_f = \dist$.
The conclusion now follows from
Corollary~\ref{c:us}(iv).
%
\end{proof}

\begin{Remark}
If $f$ is a continuous function such that $f>0$ everywhere in $\Sigma$,
by Theorem~\ref{t:unif} $\dist$ is the unique minimizer of $J$.
In this case we have
\begin{equation}\label{f:psv}
\vf>0\ \textrm{in}\ \overline{\Omega}\setminus\overline{\Sigma},
\end{equation}
and the uniqueness result also follows from Theorem 6.1
in \cite{CMg}. Nevertheless, condition (\ref{f:psv}) is only necessary
for the uniqueness.
As an example, let $\Omega\subset\R^2$ be an ellipsis
centered at the origin and with a focus in $(1,0)$.
Let $f\in C(\Omega)$ satisfy $f>0$ in
$\Omega \cap \left((-1,1)\times \R^+\right)$, $f=0$ otherwise in $\Omega$.
Let $h$ be
the indicator function of the ball $\overline{B}_1(0)$, so that $\Lambda_h=+\infty$.
In this case $\Sigma$ is the interval jointing $(-1,0)$ with $(1,0)$,
without endpoints,
so that $\Sigma \subseteq \spt(f)$, and hence, by
Theorem \ref{t:unif}, $\dist$ is the unique minimizer of $J$.
On the other hand it can be easily checked that in this case $\vf>0$
only in $\Omega \cap \left((-1,1)\times \R^+\right)$.
\end{Remark}

\section{Extensions and applications}
\label{s:ext}

\subsection{Some extension}

The existence result can be generalized without any effort to
minimum problems of the form
\begin{equation}\label{f:Jg}
\min_{u\in\Wuu(\Omega)}J(u)\,, \qquad
J(u) := \int_{\Omega} [h(Du) - g(x,u)]\,
dx\,,
\end{equation}
where
\begin{itemize}
\item[(H4)]
$g\colon\Omega\times\R\to\R$ is a measurable
function, Lipschitz continuous with respect to
the second variable, satisfying
$g(\cdot, 0)\in L^1(\Omega)$ and
\[
0\leq D_u g(x,u)\leq L,\qquad
\textrm{a.e.}\ (x,u)\in\Omega\times\R\,.
\]
\end{itemize}

\begin{Theorem}
\label{t:exig}
Let $\Omega\subset\R^n$ be a set fulfilling (\ref{f:Omega}).
Assume that (H1), (H2) and (H4) hold, together with
\begin{itemize}
\item[(H3')]
there exists
$H^*\leq H_0$
such that
\begin{equation}\label{f:cres2}
L\, c(H^*,\inradius) \leq \Lambda.
\end{equation}
\end{itemize}
Then
the distance function $\dist$
is a minimizer of $J$ in $W^{1,1}_0(\Omega)$.
Moreover, if $u\in W^{1,1}_0(\Omega)$ is a
minimizer of $J$ in $W^{1,1}_0(\Omega)$, then
$Du(x)\in K$ for a.e.\ $x\in\Omega$.
\end{Theorem}

\begin{proof}
See the proof Theorem~5.3 in \cite{CCCG}.
\end{proof}

\begin{Remark}
If $g(x,\cdot)$ is concave for a.e.\ $x\in\Omega$,
assumption (H4) can be relaxed with the following
requirement:
\begin{equation}\label{f:gconvex}
g'_+(x,0)\leq L,\quad
g'_-(x,\dist(x)) \geq 0,
\qquad \textrm{a.e.}\ x\in\Omega.
\end{equation}
Namely, let us define
\[
\tilde{g}(x,u) :=
\begin{cases}
g(x,u),
&\textrm{if $0 \leq u \leq \dist(x)$},\\
g(x,0)+g'_+(x,u)\,u,
&\textrm{if $u<0$},\\
g(x,\dist(x))+g'_-(x,\dist(x))\,(u-\dist(x)),
&\textrm{if $u>\dist(x)$}.
\end{cases}
\]
Let $\widetilde{J}$ denote the functional in (\ref{f:Jg})
with $\tilde{g}$ instead of $g$.
It is plain that, if $g$ satisfies (\ref{f:gconvex}),
then $\tilde{g}$ satisfies (H4),
so that by Theorem~\ref{t:exig} $\dist$ is a minimizer of $\widetilde{J}$.
Moreover, $\tilde{g}\geq g$, so that $\widetilde{J}\leq J$,
and $J(\dist) = \widetilde{J}(\dist)$,
hence $\dist$ is a minimizer of $J$.
\end{Remark}

Concerning the uniqueness of the minimizer, the result is the
following.

\begin{Theorem}\label{t:unig}
Assume that (H1), (H2), (H3') and (H4) hold, and let
\[
A := \{x\in\Omega;\
g(x,\dist(x)) > g(x,u)\
\forall u < \dist(x)\}\,.
\]
If $\Sigma\subset\overline{A}$, then $\dist$ is the
unique minimizer of $J$.
\end{Theorem}

\begin{proof}
Let $u_A$ be the function defined by (\ref{f:maxz})
with $S=\overline{A}$.
Let $u\in\Liprz$ be a minimizer of $J$.
It is clear that we must have $u=\dist$ in $A$,
for otherwise we would have $J(u) > J(\dist)$.
{}From Proposition~\ref{l:us}(iii)
we have that $u\geq u_A$.
Moreover, from Proposition~\ref{l:us}(v) we conclude that
$\dist = u_A\leq u\leq\dist$,
that is $u=\dist$.
\end{proof}

\begin{Remark}
The analogous of Theorems~\ref{t:exig} and \ref{t:unig}
can be proved if in (H4) we require that
$-L \leq D_u g(x,u)\leq 0$
for a.e.\ $(x,u)\in\Omega\times\R$.
\end{Remark}

\subsection{A system of PDEs of Monge-Kantorovich type}
\label{s:MK}
Let us read the results of Sections~\ref{s:exis}
and~\ref{s:uni} in terms of
properties of solutions to
the Monge-Kantorovich system of PDEs
(\ref{f:syst1})-(\ref{f:syst2}).

Assume that $f\colon\Omega\to\R$ is a bounded continuous non negative
function, and that $K$ satisfies (\ref{f:ipoK}).
Let $\vf$ be the function defined in (\ref{f:vf}).
Then $\vf$ is continuous (see \cite{CMf}, Theorem~7.2)
and $\{\vf>0\}\subseteq\Tf$ (see (\ref{f:vanv})).

We claim that $u$ is a solution to (\ref{f:minfak}) if and only if
there exists a bounded function $v\in C(\Omega)$ such that
$(u,v)$ is a solution to (\ref{f:syst1})-(\ref{f:syst2}).
Namely, let $u$ be a solution to (\ref{f:minfak}).
Then, by Theorem~\ref{t:uf}
and Corollary~\ref{c:us}, we have that
$u_f\leq u\leq\dist$, $u=\dist$ on the set $\{\vf>0\}$,
and $Du = D\dist$ a.e.\ in $\{\vf>0\}$.
Since $(\dist, \vf)$ is a solution to (\ref{f:syst1})-(\ref{f:syst2})
(see \cite{CMf}, Theorem~7.2),
then the above properties imply that also $(u, \vf)$
is a solution of the same system.
Conversely,
if $(u,v)\in \Lip(\Omega) \times C(\Omega)$, with $v$ bounded, is a solution
to (\ref{f:syst1})--(\ref{f:syst2}), then
$u$ is a solution to~(\ref{f:minfak})
(see (\ref{f:estih})),
and $v=\vf$ (see \cite{CMg}, Theorem~6.1).

Finally, Theorem~\ref{t:unif} states that the system (\ref{f:syst1})-(\ref{f:syst2})
admits the unique solution $(\dist, \vf)$
if and only if $\Sigma \subset \spt(f)$.

\subsection{Convergence of solutions to $p$-Laplace equation}

Let $f\in L^{\infty}(\Omega)$, $f\geq 0$.
In \cite{CMh} it is proved that the functionals
\[
J_p(u) :=
\begin{cases}
\displaystyle\int_\Omega \left[\frac{1}{p}\,\gauge(Du)^p-fu\right]\, dx & u\in W^{1,p}_0(\Omega), \\
+\infty & \textrm{otherwise in}\ L^1(\Omega)
\end{cases}
\]
$\Gamma$--converge in $L^1(\Omega)$ as $p\to \infty$ to the functional
\[
J(u) :=
\begin{cases}
-\int_\Omega fu \, dx\,, & u\in \Liprz, \\
+\infty & \textrm{otherwise in}\ L^1(\Omega)
\end{cases}
\]
and that the sequence $(u_p)$ of the minimizers of $J_p$ is
bounded in $W^{1,q}_0(\Omega)$ for every $q>1$.
As a consequence, any converging subsequence of $(u_p)$ converges
to a minimizer of $J$.
We remark that, for every $p>1$, $u_p$ is the unique distributional solution
of the anisotropic $p$-Laplace equation
\[
\begin{cases}
-\dive(A_p(Du_p)) = f
&\textrm{in $\Omega$},\\
u_p = 0
&\textrm{on $\partial\Omega$,}
\end{cases}
\]
where $A_p(0)=0$ and $A_p(\xi) = \gauge(\xi)^{p-1}D\gauge(\xi)$
for every $\xi\neq 0$.
By Theorem~\ref{t:unif}, if $\Sigma \subset \spt(f)$, then
$\dist$ is the unique minimizer of $J$.
In this case we can conclude that
the whole sequence $(u_p)$ converges weakly in $W^{1,q}_0(\Omega)$
to $\dist$.

In the special case $\gauge(\xi)=|\xi|$, the minimizer
$u_p$ of $J_p$ is the unique distributional solution in
$W^{1,p}_0(\Omega)$ of the $p$-Laplace equation
\[
\begin{cases}
-\Delta_p u = f
&\textrm{in $\Omega$},\\
u = 0
&\textrm{on $\partial\Omega$}.
\end{cases}
\]
In this case it is known (see \cite{BDM})
that $(u_p)$ converges in $C(\Omega)$ to the function
\[
U(x) :=
\begin{cases}
\dist(x),
&\textrm{if $x\in\overline{\Omega}\setminus A$},\\
w(x),
&\textrm{if $x\in A$,}
\end{cases}
\]
where $A := \inte\{x\in\Omega;\ f(x) = 0\}$, and
$w\in C(\overline{A})$ is the unique viscosity solution
of the $\infty$-Laplace equation
\[
\begin{cases}
-\Delta_{\infty} w = 0
&\textrm{in $A$},\\
w = \dist
&\textrm{on $\partial A$}
\end{cases}
\]
(see \cite{Jen,BaBu}).
If $\Sigma\subset\spt(f)$,
then our results state that, for every $q>1$, the sequence
$(u_p)$ converges to $\dist$ as $p\to\infty$
in $W^{1,q}$.
As a consequence, we have $w=\dist$ in $A$.


\end{document}